# *Graphs of real functions with pathological behaviors*
Claudio Bernardi
*Sapienza*, University of Rome

*A friend of mine* -  This paper is dedicated to the memory of Franco Montagna.

Franco was a dear friend of mine. I guess that, among the people who took part in the *Conference in memoriam Franco Montagna*, I am the person who has known him for the longest time: we were students together at the university of Pavia in the late Sixties. For our "*tesi di laurea*", in the early Seventies we went to Ferrara together to speak with Roberto Magari (at that time he was at the University of Ferrara), who was the thesis advisor for both of us.

One year later Roberto Magari moved to Siena, and Franco and I obtained a fellowship from the *CNR*, the Italian science foundation. For many years Franco and I were colleagues and friends in Siena. I wish to recall our first trip to the United States: that time, as well, the two of us went to New York together for a couple of weeks.

When I was in Siena, I was invited to have dinner at Franco and Tonina's house many times, I would say uncountably many times (I actually could never count the number of times I went there). In those years Tonina and Franco were true friends for me. Today I would like to say:

*thank you, Franco, wherever you are,*
*thank you, Tonina.*

Then I moved to the University of Rome. Over the past twenty-five years unfortunately I only met Franco occasionally. The last time was in Torino in the summer of 2014, when I had the opportunity to speak with Franco a lot. A few months later I was told about his illness.

**1.** *Unusual periodic functions* -  I intend to present and discuss functions from **R** to **R** whose graphs have unexpected properties. So the subject is not mathematical logic. However, I think that there are connections with logic: I will be looking for counterexamples, paradoxes, pathological objects, analogies, and this always gives, in some sense, a logical flavor to the subject. In general, I will discuss situations where intuition, and also perception, may conflict with a rigorous understanding of what is going on.

Let me start with an elementary question.
Let *f* be a *non-constant periodic function*. Is it always true that *f* admits a minimum positive period?
If we think about a usual trigonometric function, like sin (*x*), the answer is *yes* (and indeed we often speak of *the* period of the function). However, it is enough to think about the Dirichlet function to find out that any rational number (different from 0) is a period.

To construct other examples of a periodic function without a minimum period, first I restrict the domain of functions: instead of **R** I will consider the algebraic extension **Q**[$\sqrt{2}$] = {$a + b\sqrt{2} \mid a, b \in$ **Q**} $\subseteq$ **R**; later I will consider functions that have domain **R**, but the Axiom of Choice will be needed.

The point is that, on the one hand, the additive group **Q**[$\sqrt{2}$] is obviously isomorphic to **Q**×**Q**; but, on the other hand, **Q**[$\sqrt{2}$] is an ordered subset of **R**: this also is obvious, but the order relation in **Q**[$\sqrt{2}$] has nothing to do with the above isomorphism.

Let us define two functions $p$ and $q$ from **Q**[$\sqrt{2}$] into itself:

$$p(a + b\sqrt{2}) = a \qquad \text{and} \qquad q(a + b\sqrt{2}) = b\sqrt{2}.$$

From an algebraic point of view, these functions are very simple: in particular, if **Q**[$\sqrt{2}$] is seen as a vector space over **Q**, $p$ and $q$ are linear maps; and, if **Q**[$\sqrt{2}$] is regarded as an additive group, so as isomorphic to **Q**×**Q** (more precisely, to **Q**×**Q**[$\sqrt{2}$]), they are simply the two projections.

Now, given any rational number $r$, it is very easy to verify that $p(a + b\sqrt{2} + r\sqrt{2}) = a = p(a + b\sqrt{2})$ (adding $r\sqrt{2}$ the value of $p$ does not change), and similarly $q(a + b\sqrt{2} + r) = b\sqrt{2} = q(a + b\sqrt{2})$. So, $p$ and $q$ are periodic, without a minimum positive period: any rational multiple of $\sqrt{2}$ is a period for $p$ and any rational number is a period for $q$.

Of course, the graph of $p$, as well as the graph of $q$, repeats itself just like the graph of any periodic function; however, from an intuitive point of view the situation may not be clear, also because we cannot represent the graphs of $p$ and $q$ in the usual way. Indeed, as $a$ or $b$ or both increase, the number $x = a + b\sqrt{2}$ also increases. On the other hand, as $x = a + b\sqrt{2}$ increases, one might expect that also $a$ and $b\sqrt{2}$ tend to increase. Occasionally, it can happen that $x$ increases and $a$ (or $b\sqrt{2}$) decreases, but, if this is the case, then the increase in $b\sqrt{2}$ (resp. $a$) must be greater than the decrease in $a$ (resp. $b\sqrt{2}$). However this is not completely correct: in fact, the functions $p$ and $q$ are periodic, so they repeat themselves. This means that for neither of them it makes sense to say that they "tend to increase".

Moreover we have:

$$p + q \text{ is the identity function, since } (p + q)(a + b\sqrt{2}) = a + b\sqrt{2}.$$

The fact that the identity can be expressed as the sum of two periodic functions was pointed out by Mortola and Peirone in [9] (they referred to **R** and used the Axiom of Choice).

**2. *Functions whose graphs repeat themselves at different levels*** - Now, let me consider a function like $g(x) = \sin(x) + x/2$ (see Figure 1; I considered $x/2$ instead of $x$ only because the graph of $\sin(x) + x$ is less clear).

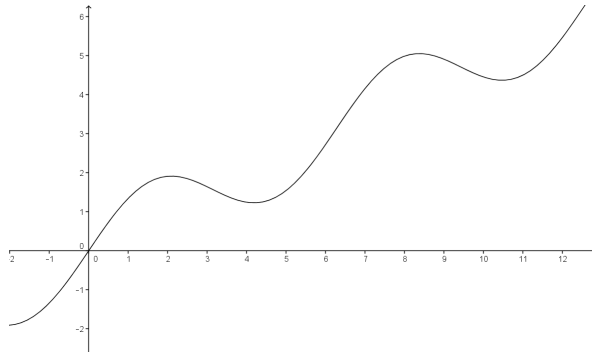

Figure 1

The function *g* is not periodic, it is said to be *quasiperiodic*: that means that there exist two numbers *T* and *c*, both different from 0, such that $g(x + T) = g(x) + c$ for every *x*. *T* is called a *quasiperiod*.

The graph of a quasi periodic function also repeats itself, not at the same "level", but each time a little higher or a little lower. More precisely, if *T* and *c* have the same sign, then the graph will repeat itself a little higher, while if *T* and *c* have different signs, then the graph will repeat itself a little lower; one could speak respectively of an *increasing* or a *decreasing* quasiperiodic function. For instance, the function $h(x) = \sin(x) - x/2$ (see Figure 2) is a decreasing quasiperiodic function.

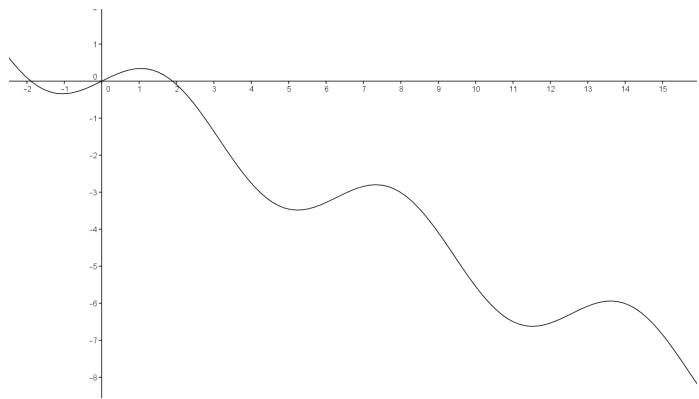

Figure 2

Contrasting the intuitive meaning of a quasiperiodic function, we have that *the functions p and q above are both not only periodic but also, at the same time, increasing quasiperiodic, and decreasing quasiperiodic*. For instance:

$$p(a + b\sqrt{2} + 1) = a + 1 \quad \text{(adding 1 the value of } p \text{ increases by 1);}$$
$$p(a + b\sqrt{2} + \sqrt{2} - 1) = a - 1 \quad \text{(adding } \sqrt{2} - 1 \text{ the value of } p \text{ decreases by 1).}$$

Similar remarks apply to *q*. More precisely, *any* non-zero number belonging to $\mathbf{Q}[\sqrt{2}]$ is a period or a quasiperiod, both for *p* and for *q*.

Notice also that the graphs of *p* and *q* are both *dense* in $\mathbf{R}^2$. This is not hard to be proved.

**3. *Additive functions*** -  In order to find functions with similar properties, but whose domain is the whole set **R**, we start with the functional equation

$$f(x + y) = f(x) + f(y) \qquad \text{for every } x, y.$$

A function is said to be *additive* if it satisfies this equation. It is easy to see that a function is additive *iff* it is a linear map when **R** is seen as a vector space over **Q**.
Cauchy proved that every continuous additive function is of the kind $f(x) = kx$, for some value of the parameter $k$. Of course, if this is the case, the graph of $f$ must be a straight line. About a century ago, using the Axiom of Choice Hamel found *discontinuous* additive functions with domain **R**. An easy construction is as follows. Applying Zorn's Lemma we can find two subspaces **G** and **H** of **R** (where **R** is seen again as a vector space over **Q**), such that **R** = **G** ⊕ **H** (where the symbol ⊕ denotes direct sum). This means that each real number $x$ can be uniquely written as a sum $x = g + h$ with $g \in$ **G** and $h \in$ **H**. So, again, we can define the projections $p(x) = g$ and $q(x) = h$. The domain of $p$ and $q$ is **R**; they enjoy the same pathological properties we have seen before, namely:
-  $p$ and $q$ are periodic and, at the same time, quasiperiodic;
-  every number (different from 0) is a period or a quasiperiod;
-  $p + q$ is the identity function;
-  the graphs of $p$ and of $q$ are dense in the plane.

**4. *Homogeneous graphs*** - Let us go back to Euclid's *Elements*. The celebrated definition 4 in the first Book of *Elements* is as follows:

*"a straight line is a line which lies evenly with the points on itself."*

The precise meaning of Euclid's statement has been widely discussed. The usual interpretation is explained by Heath, who says [5 - p. 167]: *"we can safely say that the sort of idea which Euclid wished to express was that of a line which presents the same shape at and relatively to all points on it, without any irregular or unsymmetrical feature distinguishing one part or side of it from another."* Hence I simply intend Euclid's definition as: in a straight line no point is privileged, or, in more technical words, a straight line is *homogeneous*.
I claim that the same property holds for the graph of any additive function $f$:

*for every A, B belonging to the graph, there is an isometry φ (in fact a translation) such that φ maps the graph onto itself and φ(A) = B.*

The proof is simple. Let $A = (x_A, f(x_A))$ and $B = (x_B, f(x_B))$ be two points of the graph of $f$. The translation by the vector $AB$ maps a point $(x, f(x))$ of the graph into the point $(x + x_B - x_A, f(x) + f(x_B) - f(x_A)) = (x + x_B - x_A, f(x + x_B - x_A))$, which is in turn a point of the graph.
We can conclude that all points on the graph of any additive function look the same, in the sense that any two points cannot be distinguished from each other within the graph. In particular, even if the points of the graph have *y*-values which are "very far from each other", no point is in a higher or lower position relative to the other points.

I wish to stress an interesting point: discontinuous additive functions present a very irregular behavior (they are nowhere continuous, unbounded in any interval, ...) and, at the same time, they are regular like a straight line.

**5.** *Everywhere surjective functions* - Let me introduce another class of functions. A function from **R** to **R** is said to be an *everywhere surjective function* (or an *everywhere surjection*, or also a *strong Darboux function*) if, for every interval (*a*, *b*) and for every *y*, there is an $x \in (a, b)$ such that $f(x) = y$.

In other words, a function from **R** to **R** is an everywhere surjection if its restriction to *any* interval is surjective. Notice that an equivalent definition is the following: $f^{-1}\{y_0\}$ is dense in **R** for every $y_0$.

At first sight, one might doubt that such a function exists.

The first example of an everywhere surjection was introduced by Lebesgue in 1904 ([8] - p. 97). Recently, this concept has been generalized and deeply studied (see for instance [2] and [7]).

It is interesting to recall why Lebesgue introduced these functions. The *intermediate value theorem* states that, if *f* is a continuous function defined in an interval [*a*, *b*], then *f* takes on any value between *f*(*a*) and *f*(*b*) at some point within the interval. Sometimes this property has been confused with the definition of a continuous function: some ancient books gave the following definition of continuity: «a function *f* from **R** to **R** is continuous if, in any interval [*a*, *b*], *f* takes on any value between *f*(*a*) and *f*(*b*)». This "definition" was criticized by Darboux in a famous paper in 1875. And Lebesgue said: «*On me permettra de signaler qu'en 1903 on enseignait encore dans un lycée de Paris la définition critiquée dès 1875 par Darboux. Cela est d'autant plus étonnant que la propriété qui est énoncée dans la définition de Cauchy est celle qui intervient directement dans presque toutes les démonstrations.*»

To show that the property stated in the conclusion of the intermediate value theorem is not equivalent to continuity, Lebesgue gave an example of what we now call an everywhere surjection, which obviously in any interval takes on any value between $f(a)$ and $f(b)$ without being continuous. In fact, an everywhere surjection cannot be continuous at any point $x_0$: it is enough to observe that, if a function *f* is continuous at $x_0$, then there are two positive number $\varepsilon$, $\delta$ such that $f(x_0 - \delta, x_0 + \delta) \subseteq (f(x_0) - \varepsilon, f(x_0) + \varepsilon)$; but, if *f* is an everywhere surjection, then $f(x_0 - \delta, x_0 + \delta) = \mathbf{R}$. In the next Section we will see two constructions of everywhere surjective functions; the second of these constructions is not too different from Lebesgue's one.

Let us discuss a little more the previous quotation and compare the usual correct definition of continuity with the one criticized by Darboux and Lebesgue. Of course Lebesgue was right: the "definition" according to which a function *f* is continuous if, in any interval [*a*, *b*], it takes on any value between *f*(*a*) and *f*(*b*) within the interval, is wrong, in the sense that it does not capture what we want. Moreover, as Lebesgue observed, if we intend to use the definition in a deductive mathematical framework, only the definition of continuity that we know to be correct works. On the other hand, if our aim were *only* to provide a good

*description* of the behavior of a graph, the "wrong definition" is simpler and more intuitive for explaining what is going on in the graph, and in many simple cases it works.

**6. *Constructions of everywhere surjective functions* -** Examples of an everywhere surjection can be given following different procedures. For instance, one can refer to the Cantor set and proceed as follows (see [6, Example 2.2]).

Let $\{(a_i, b_i) \mid i \in \mathbf{N}\}$ be a denumerable open basis for $\mathbf{R}$. We consider copies $C_i$ (where $i \in \mathbf{N}$) of the *ternary Cantor set*, chosen inductively so that $C_i \subset (a_i, b_i)$ for all $i$ and every $C_i$ is disjoint from $\cup_{j<i} C_j$ (at any step there is enough room because the Cantor set has Peano-Jordan measure zero). Let $f_i$ be a function from $C_i$ onto $\mathbf{R}$. Finally, we define $f(x) = f_i(x)$ if $x \in C_i$, and $f(x) = 0$ otherwise. Clearly $f$ is an everywhere surjection, since any interval $(a, b)$ contains a set $C_i$ for some $i$ (the Axiom of Choice is not needed).

Also John Conway, well-known for his original and creative approach to many mathematical fields, proposed an example of an everywhere surjection: *the base* 13 *function* - see [11] and also [10]. I will not examine here Conway's procedure to define an everywhere surjection $h$, but I present a different construction, which is similar but simpler. We can write real numbers using *ternary notation*. Let $q$ be integer part of a number $x$; we pay attention only to the last two digits "2" in the ternary expansion of $x$ after the decimal point (if there is at most one digit "2" or there are infinitely many digits "2", then $h(x) = 0$). So

$$x = q. \ldots \underbrace{2b_1b_2\ldots b_m}\underbrace{2y_1y_2}\ldots$$

whereas all digits $b_1b_2\ldots b_m$ and $y_1y_2\ldots$ are either 0 or 1. Then $h(x)$ is defined to be the number having $b_1\ldots b_m$ as its integer part and $y_1y_2\ldots$ as its fractional part, where both these expressions are intended in *binary notation*:

$$h(x) = b_1b_2\ldots b_m \bullet y_1y_2\ldots$$

It is not hard to check that $h$ is an everywhere surjection; notice that the first digits of $x$ after the decimal point (represented by the dots "..." before the two digits "2") guarantee that a given value is assumed in any given interval.

To obtain also negative numbers in the image of $h$, we can refer to the digit $b_1$: for instance, we can decide that $f(x)$ is negative if $b_1$ is 0, and positive otherwise.

Also notice that this function $h$
- maps rational numbers into rational numbers,
- it is periodic, with period 1.

**7. *Some results* -** Let me briefly state some theorems about additive functions and everywhere surjective functions (see for instance [3]). Some of them are simple, while others require a not completely obvious proof.

- The graph of an everywhere surjective function is dense in $\mathbf{R}^2$;

- an additive function is periodic *iff* it is not injective;
- an additive function is an everywhere surjection *iff* it is surjective but not injective;
- the graph of an additive function is symmetric with respect to any of its points, like a straight line.

More technical results show that everywhere surjections are not isolated phenomena. In particular, in [1] it is proved that there is a vector space Δ of functions from **R** to **R**, such that every non-zero element of Δ is an everywhere surjection and the dimension of Δ is 2 raised to the cardinality of the continuum.

There are also generalizations of the concept of an everywhere surjective function, which are obtained considering other sets instead of **R**. In particular, if $X$ is a topological space, a function $f : X \to X$ is said to be *everywhere surjective* if, for every non-empty open set $A$, we have $f(A) = X$.

Now we have the following characterization (see [10]). Let $X$ be a topological space of cardinality $k$ (where $k$ is infinite); there exists an everywhere surjective function from $X$ to $X$ *iff* every nonempty open subset of $X$ has cardinality $k$ and there are $k$ subsets of $X$, which are dense in $X$ and pairwise disjoint.

Moreover, considering functions from **R** to **Q**, one can find an example of an everywhere surjective function such that, for any $y \in \mathbf{Q}$, the inverse image of $y$ has a *positive* measure (C. Rainaldi, unpublished).

We conclude with a topological remark. In this context, we have simple examples of connected but not pathwise connected plane sets (in fact, with uncountably many path-components). Indeed, consider the region below the graph of a function, that is, the set $H = \{(x, y) \mid y \leq f(x)\}$; sometimes $H$ is called *hypograph* of $f$. Notice there are functions with domain **R** whose hypograph is not connected: it is enough to think of the function defined as follows: $f(x) = 1/x$ if $x$ is positive and $f(x) = 0$ otherwise.

Now, if $f$ is either a discontinuous additive function or an everywhere surjection, we have that (see [4]):
- the hypograph of $f$ is *connected* (this statement is intuitive, but the proof is not trivial); more generally, if the graph of a function is dense in the plane, then the hypograph is connected;
- the hypograph of $f$ is *not pathwise connected*: indeed the hypograph of $f$ has uncountably many path-connected components, which are all the half-lines parallel to the $y$-axis.

**Some References**